\documentclass[reqno]{amsart}
\usepackage{times}
\usepackage{psfrag}
\usepackage[T1]{fontenc}
\usepackage{amssymb}

\makeatletter
 \theoremstyle{plain}    
 \newtheorem{thm}{Theorem}[section]
 \numberwithin{equation}{section} %% Comment out for sequentially-numbered
 \numberwithin{figure}{section} %% Comment out for sequentially-numbered
 \theoremstyle{plain}
 \theoremstyle{remark}
 \newtheorem{rem}[thm]{Remark}
 \theoremstyle{definition}
  \newtheorem{example}[thm]{Example}
 \theoremstyle{definition}
 \newtheorem{defn}[thm]{Definition}
 \theoremstyle{plain}    
 \newtheorem{prop}[thm]{Proposition} %%Delete [thm] to re-start numbering
 \theoremstyle{plain}    
 \newtheorem{lem}[thm]{Lemma} %%Delete [thm] to re-start numbering

 \usepackage{bbm} 
\usepackage{graphicx}
\newcommand{\Q}{\mathbb {Q}}

\newcommand{\1}{\boldsymbol{1}}

\newcommand{\N}{\mathbb {N}}

\newcommand{\I}{\mathbb {I}}

\AtBeginDocument{

}
\makeatother
\begin{document}
\let\language\relax
\title[A distributional limit law for the continued fraction digit sum]{A distributional limit law for the\\ continued fraction digit sum} 
\author{Marc Kesseböhmer and Mehdi Slassi}
\email{mhk@math.uni-bremen.de, slassi@math.uni-bremen.de}
\date{15 June, 2005}
\begin{abstract}
We consider the continued fraction digits as random variables measured
with respect to Lebesgue measure. The logarithmically scaled and normalized
fluctuation process of the digit sums converges strongly distributional
to a random  variable uniformly distributed on the unit interval.
For this process normalized linearly we determine a large deviation
asymptotic. 
\end{abstract}

\address{Fachbereich 3 - Mathematik und Informatik, Universität Bremen, Bibliothekstrasse
1, 28359 Bremen, Germany. }
\maketitle

\section{Introduction and statement of main results}\let\language\relax

Any number $x\in\I:=\left[0,1\right]\setminus\Q$ has a simple infinite
continued fraction expansion\[
x=\frac{1}{a_{1}\left(x\right)+{\displaystyle \frac{1}{a_{2}\left(x\right)+\cdots}}},\]
where the unique \emph{continued fraction digits} $a_{n}\left(x\right)$
are from the positive integers $\mathbb{N}$. The \emph{Gauss transformation
$G:\I\rightarrow\I$} is given by\[
G(x):=\frac{1}{x}-\left\lfloor \frac{1}{x}\right\rfloor ,\]
where $\left\lfloor x\right\rfloor $ denotes the greatest integer
not exceeding $x\in\mathbb{R}$. Write $G^{n}$ for the $n$-th iterate
of $G$, $n\in\mathbb{N}_{0}=\left\{ 0,1,2,\ldots\right\} $ with
$G^{0}=\textrm{id.}$ It is then well known that for all $n\in\mathbb{N}$,
we have \[
a_{n}(x)=\left\lfloor \frac{1}{G^{n-1}x}\right\rfloor .\]
 Clearly, the $a_{n}$, $n\in\mathbb{N}$, define random variables
on the measure space $\left(\I,\mathcal{B},\mathbb{P}\right)$, where
$\mathcal{B}$ denotes the Borel $\sigma$-algebra of $\I$ and $\mathbb{P}$
some probability measure on $\mathcal{B}$. Then each $a_{n}$ has
infinite expectation with respect to the Lebesgue measure on $[0,1]$,
which we will denote by $\lambda$. By the ergodicity of the Gauss
transformation with respect to the Gauss measure $d\mu\left(x\right):=\frac{1}{\log2}\frac{1}{1+x}d\lambda\left(x\right)$
we readily reproduce Khinchin's result on the geometric mean of the
continued fraction digits, i.e.\[
\sqrt[n]{a_{1}\cdots a_{n}}\to K,\,\lambda\textrm{-a.e.,}\]
where $K=2,685...$ denotes the Khinchin constant. A similar result
holds for the harmonic mean. Also by a classical result of Khinchin
(cf. \cite{Khinchin:36}), we know that for $\lambda$-almost every
$x\in\left[0,1\right]$ we have for infinitely many $n\in\N$ \[
a_{n}\left(x\right)>n\log n.\]
From this Khinchin deduced that the arithmetic mean of the continued
fraction digits is divergent, i.e. \[
\lim_{n\to\infty}\frac{S_n}{n}
=\infty\quad\;\lambda\textrm{-a.e.},\] where $S_n(x) :=a_{1}(x)+\cdots+a_{n}
\left(x\right)$, $x\in \I$.
It was again Khinchin who showed that nevertheless for a suitable 
normalising sequence a weak law of large numbers holds. That is
$\frac{S_n}{(n\log n)}$ converges in measure to $1/\log 2$ with respect
to $\lambda$. However, according to Phillip \cite{Philipp:88} there is no (reasonable)
normalising sequence $(n_k)$ with $(n_k/k)$ non-decreasing such that a strong law of 
large numbers is satisfied. More precisely, we either have
\[
\sum _{k=1}^\infty \frac{1}{n_k}<\infty \;\mathrm{ and }\; \lim_{k\to \infty}
\frac{S_k}{n_k}=0\; \;\lambda\textrm{-a.e.},\;\;\;\mathrm{ or }\;\;\;
\sum _{k=1}^\infty \frac{1}{n_k}=\infty \;\mathrm{ and }\; \lim_{k\to \infty}
\frac{S_k}{n_k}=\infty \;\; \lambda \textrm{-a.e.}
\]
In contrast to this Diamond and Vaaler showed in \cite{DiamondVaaler:86}
that for the trimmed sum $S^\flat_n := S_n - \displaystyle{\max_{1\leq \ell\leq n}} a_\ell$
a strong law of large numbers holds in the sense that   
\[
\lim_{n\to\infty} \frac{S^\flat_n}{n\log n}=\frac{1}{\log 2}\;\;\lambda\textrm{-a.e.}
\]
This shows that the intricate stochastic properties of $S_n$ arise from the occurrences of rare
but exceptionally large continued fraction digits.  
L{\'e}vy in  \cite{Levy:52} was the first to predict  non-degenerated limit laws in the context of continued
fractions  -- namely stable laws. Actually, we have that $(a_n)$ belongs to the domain
 of attraction to a stable law with characteristic exponent $1$. More precisely 
we have the following convergence in distribution with respect to any absolutely continuous 
probability measure $\mu\ll \lambda$  
 \[
\frac{S_k}{k/\log 2}-\log k \stackrel{\mu}
{\longrightarrow} F,\] 
where $F$ has a stable distribution.

This was first proved by Heinrich in \cite{Heinrich:87} under the assumption that $\mu$  has strictly  positive density. 
He also showed that the error term of this convergence is of order $\mathcal{O}((\log k)^2/k)$. The convergence
was reproved in \cite{Hendley:00} for arbitrary densities but without 
providing error bounds. In \cite{Philipp:88}  where also the stabel law for the sum $S_n$ is proven,   
a central limit theorem for the sufficiently trimmed sum is obtained. 

Instead of taking a global centralising sequence as in the above situation we will focus on
the pointwise behaviour of the fluctuation of $S_n$. For this we define 
\begin{equation}
X_{n}\left(x\right):=\max\left\{S_k\left(x\right)\;:\;S_k \left(x\right)\leq n,\;
 k\in\mathbb{N}_{0}\right\} ,\quad x\in\I,\label{process}
\end{equation}
and investigate the process $n-X_{n}$ of pointwise fluctuations. Employing  
infinite ergodic theory we are able to derive  non-degenerated results 
describing new aspects of the stochastic structure of $S_n$. 
The underlying dynamical system will be given by the Farey
map and the connection to metrical number theory is established via a certain 
return time process for this system. This process turns out to be conceptually
related to renewal theory which guided us to find the following main
results of this paper.

\begin{itemize}
\item \textbf{Uniform Law:} For $\boldsymbol{U}$ uniformly distributed
random variable on $\left[0,1\right]$ we have the following convergence
in distribution with respect to $\lambda$\[
{\displaystyle \frac{\log\left(n-X_{n}\right)}{\log\left(n\right)}\stackrel{\lambda}
{\longrightarrow}\boldsymbol{U}}.\]
\item \textbf{Large Deviation Law:} For all $\varepsilon\in\left(0,1\right)$
we have \[
{\displaystyle \lambda\left(\frac{n-X_{n}}{n}>\varepsilon\right)\sim\frac{-\log\varepsilon}
{\log n}},\]
where $a_{n}\sim b_{n}$ means that $\lim_{n\to\infty}\frac{a_{n}}{b_{n}}=1$.
\end{itemize}
These statements are implied by the slightly more general Theorem
\ref{thm:main} stated and proved in Subsection \ref{sub:Distributional-limit-laws}. 
Note that the \emph{uniform law} gives the following convergence in
probability\begin{equation}
{\displaystyle n-X_{n}}\stackrel{\lambda}{\to}\infty,\label{eq:infty}\end{equation}
whereas the \emph{large deviation} law in particular means \begin{equation}
{\displaystyle \frac{n-X_{n}}{n}}\stackrel{\lambda}{\to}0.\label{eq:Negelg}\end{equation}
For further interesting result in the context of continued fraction
digit sums we would like to refer the reader to 
\cite{GuivarchLeJan:93,GuevarchLeJan:96}, wherein  alternating sums of continued fraction digits
are considered.

The paper consists of two major parts. In the first part (Section
\ref{sub:Infinite-ergodic-theory}) we recall and develop some expects
from infinite ergodic theory allowing us to apply a theorem from \cite{KesseboehmerSlassi:05}
in the proof of the uniform law within this number theoretical context.
In the second part we study the Farey map as an infinite measure preserving
transformation and use its connections to the fluctuation process
in question to finally give the proof of both the uniform law and
the large deviation asymptotic.

\section{Infinite ergodic theory. \label{sub:Infinite-ergodic-theory} }

In the first subsection of this section we recall some basics definitions
and facts from infinite ergodic theory and give new sufficient conditions
for a set to be uniformly returning. In the second subsection we state
the general limit law which we will apply in the context of continued
fractions. Also, as an example we consider a certain interval map.

\subsection{Preliminaries}

With $\left(X,T,\mathcal{A},\mu\right)$ we always denote a conservative
ergodic measure preserving dynamical systems where $\mu$ is an infinite
$\sigma$-finite measure. For a good overview, further definitions
and details we refer the reader to \cite{Aaronson:97}.

Let \[
\mathcal{P}_{\mu}:=\left\{ \nu:\nu\:\textrm{probability measure on}\:
\mathcal{A}\,\textrm{with }\nu\ll\mu\right\} \]
 denote the set of probability measures on $\mathcal{A}$ which are
absolutely continuous with respect to $\mu$. The measures from $\mathcal{P}_{\mu}$
represent the admissible initial distributions for the processes associated
with the iteration of $T$. The symbol $\mathcal{P}_{\mu}$ will also
be used for the set of the corresponding densities.

Note that for such dynamical systems the mean return time to sets
of finite positive measure is infinite leading to the notion of wandering
rate. For a fixed set $A\in\mathcal{A}$ with $0<\mu\left(A\right)<\infty$
we set 
\[
K_{n}:=\bigcup_{k=0}^{n}T^{-k}A\quad\mathrm{and}\quad W_{n}:=W_{n}\left(A\right):=\mu\left(K_{n}\right),\qquad n\geq0,\]
and call the sequence $\left(W_{n}\left(A\right)\right)$ the \emph{wandering
rate} of $A$. Note that for the wandering rate the following identity
holds\[
W_{n}\left(A\right)=\sum_{k=0}^{n}\mu\left(A\cap\{\varphi>k\}\right),\]
where \[
\varphi(x):=\inf\{ n\geq1:\; T^{n}(x)\in A\},\quad x\in X,\]
denotes the \emph{first return time} to the set $A$. 
Since $T$ is conservative and ergodic, for all $\nu\in\mathcal{P_{\mu}}$,\[
\lim_{n\to\infty}\nu\left(K_{n}\right)=1\quad\textrm{and}
\quad\nu\left(\left\{ \varphi<\infty\right\} \right)=1.\]
To understand the stochastic properties of a nonsingular transformation
of a $\sigma$-finite measure space one often has to study the long-term
behaviour of the iterates of its \emph{transfer operator} \[
\hat{T}:L_{1}\left(\mu\right)\longrightarrow L_{1}\left(\mu\right),\; 
f\longmapsto\hat{T}\left(f\right):=\frac{d\left(\nu_{f}\circ T^{-1}\right)}{d\mu},\]
where $\nu_{f}$ denotes the measure with density $f$ with respect
to $\mu$. Clearly, $\hat{T}$ is a positive linear operator characterised
by\[
\int_{B}\hat{T}\left(f\right)\; d\mu=\int_{T^{-1}\left(B\right)}f\; 
d\mu,\qquad f\in L_{1}\left(\mu\right),\quad B\in\mathcal{A}.\]

The ergodic properties of $(X,T,\mathcal{A},\mu)$ can be characterised
in terms of the transfer operator in the following way (cf. 
\cite[Proposition 1.3.2]{Aaronson:97}).
A system is conservative and ergodic if and only if for all 
$f\in L_{1}^{+}\left(\mu\right):=\left\{ f\in L_{1}\left(\mu\right):\; 
f\geq0\;\mathrm{and}\;\int_{X}f\; d\mu>0\right\} $
we have $\mu$-a.e. \[
\sum_{n\geq0}\hat{T}^{n}\left(f\right)=\infty.\]
Invariance of $\mu$ under $T$ means $\hat{T}\left(\1\right)=\1$.

The following two definitions are in many situation crucial within
infinite ergodic theory.
\begin{itemize}
\item A set $A\in\mathcal{A}$ with $0<\mu\left(A\right)<\infty$ is called
\emph{uniform for} $f\in P_{\mu}$ if there exists a sequence $\left(b_{n}\right)$
of positive reals such that\[
\frac{1}{b_{n}}\sum_{k=0}^{n-1}\hat{T}^{k}\left(f\right)\;\longrightarrow\;1
\qquad\mu-\textrm{a.e. \; uniformly\; on}\; A\]
(i.e. uniform convergence in $L_{\infty}\left(\mu|_{A\cap\mathcal{A}}\right)$).
\item The set $A$ is called a \emph{uniform} set \emph{}if it is uniform
for some $f\in P_{\mu}$.
\end{itemize}
The concept of regularly varying functions and sequences plays a central
role in infinite ergodic theory (see also \cite{BinghamGoldieTeugels:89}
for a comprehensive account). 

A measurable function $R:\mathbb{R}^{+}\rightarrow\mathbb{R}$ with
$R>0$ on $\left(a,\infty\right)$ for some $a>0$ is called \emph{regularly
varying} at $\infty$ with exponent $\rho\in\mathbb{R}$ if
\[
\lim_{t\to\infty}\frac{R\left(\lambda t\right)}{R\left(t\right)}=
\lambda^{\rho}\quad{\rm for\; all}\;\lambda>0.\]
A regularly varying function $L$ with exponent $\rho=0$ is called
\emph{slowly varying} at $\infty$, i.e.\[
\lim_{t\to\infty}\frac{L\left(\lambda t\right)}{L\left(t\right)}=1\quad{\rm for\; all}\;\lambda>0.\]
 Clearly, a function $R:\mathbb{R}^{+}\rightarrow\mathbb{R}$ is regularly
varying at $\infty$ with exponent $\rho$$\in\mathbb{R}$ if and
only if\[
R\left(t\right)=t^{\rho}L\left(t\right),\quad t\in\mathbb{R}^{+},\]
for some slowly varying function $L$.
A \emph{sequence} $\left(u_{n}\right)$ is \emph{regularly varying
with exponent} $\rho$ if $u_{n}=R\left(n\right),\; n\geq1$, for
some $R:\mathbb{R}^{+}\rightarrow\mathbb{R}$ regularly varying at
$\infty$ with exponent $\rho$.

\begin{rem}
\label{remark0}From \cite[Proposition 3.8.7]{Aaronson:97} we know,
that $\left(b_{n}\right)$ is regularly varying with exponent $\alpha$
if and only if $\left(W_{n}\right)$ is regularly varying with exponent
$\left(1-\alpha\right)$. In this case $\alpha$ lies in the interval
$\left[0,1\right]$ and 
\begin{equation}
b_{n}W_{n}\sim\frac{n}{\Gamma\left(1+\alpha\right)\Gamma\left(2-\alpha\right)}.
\label{eq:AsymptAaronson}
\end{equation}
In the following example we introduce interval maps giving rise to
an infinite measure preserving transformation.
\end{rem}

\begin{example}(\cite{Thaler:80}, \cite{Thaler:83})
\label{ex1}Let $\xi_{1}=\left\{ B\left(k\right):\; k\in I\right\} $
be a finite or infinite family of pairwise disjoint subintervals of
$\left[0,1\right]$ such that $\lambda\left(\bigcup_{k\in I}B\left(k\right)\right)=1$.
We consider transformations $T:\left[0,1\right]\longrightarrow\left[0,1\right]$,
satisfying the following Thaler Conditions.
\end{example}
\begin{enumerate}
\item $T_{\mid B\left(k\right)}$ is twice differentiable and $\overline{TB\left(k\right)}=\left[0,1\right]$
for all $k\in I$. 
\item There exists a non-empty finite set $J\subseteq I$ such that each
$B\left(j\right),\; j\in J$, contains a unique fixed point $x_{j}$
with $T'\left(x_{j}\right)=1$ (indifferent fixed point).
\item $\left|T'\right|\geq\rho\left(\varepsilon\right)>1$ on $\bigcup_{k\in I}B\left(k\right)\diagdown\bigcup_{j\in J}\left(x_{j}-\varepsilon,x_{j}+\varepsilon\right)$
for each $\varepsilon>0$.
\item There exists $\eta>0$ such that for all $j\in J$, $T'$ is decreasing
on $\left(x_{j}-\eta,x_{j}\right)\cap B\left(j\right)$ and increasing
on $\left(x_{j},x_{j}+\eta\right)\cap B\left(j\right)$.
\item $\frac{T''}{\left(T'\right)^{2}}$ is bounded on $\bigcup_{k\in I}B\left(k\right)$
(Adler's Condition).
\end{enumerate}
As proved in \cite{Thaler:80}, \cite{Thaler:83} $T$ is conservative
and ergodic with respect to $\lambda$, and admits an infinite $\sigma$-finite
invariant measure $\mu$ equivalent to $\lambda$. The density $\frac{d\mu}{d\lambda}$
has a version $h$ of the form\[
h\left(x\right)=h_{0}\left(x\right)\prod_{j\in J}\frac{x-x_{j}}{x-u_{j}\left(x\right)},
\qquad x\in\left[0,1\right]\diagdown\left\{ x_{j}:\; j\in J\right\} ,\]
where $u_{j}=\left(T|_{B\left(j\right)}\right)^{-1},\; j\in J$, and
$h_{0}$ is continuous and positive on $\left[0,1\right]$. From \cite{thaler:95}
we know that there exists a sequence $\left(b_{n}\right)$ of positive
numbers such that for all $f\in L_{1}\left(\mu\right)$ with $fh$
Riemann-integrable on $\left[0,1\right]$ we have\[
\frac{1}{b_{n}}\sum_{k=0}^{n-1}\hat{T}^{k}\left(f\right)\;\longrightarrow\;\int f\; d\mu,\]
uniformly on compact subsets of $\left[0,1\right]\diagdown\left\{ x_{j}:\; j\in J\right\} $.
Hence, any such subset is uniform for any $f\in\mathcal{P}_{\mu}$.
In particular, if \[
T\left(x\right)=x\mp a_{j}\left|x-x_{j}\right|^{p_{j}+1}+o\left(\left|x-x_{j}\right|^{p_{j}+1}\right)\;\left(x\to x_{j}\right)\]
 with $a_{j}>0,\; p_{j}\in\mathbb{N}\;\left(j\in J\right)$, and $p=\max\left\{ p_{j}:\; j\in J\right\} $,
then we have\[
b_{n}\sim const\cdot\left\{ \begin{array}{ll}
\frac{n}{\log\left(n\right)}, & p=1,\\
n^{1/p}, & p>1.\end{array}\right.\]

To state the crucial condition for the uniform law (see (UL) in Subsection
\ref{sub:Limit-laws.}) we need the notion of uniform returning sets
introduced in \cite{KesseboehmerSlassi:05}.

\begin{defn}
A set $A\in\mathcal{A}$ with $0<\mu\left(A\right)<\infty$ is called
\emph{uniformly returning} \emph{for} $f\in P_{\mu}$ if there exists
a positive increasing sequence $\left(b_{n}\right)$ such that \[
b_{n}\hat{T}^{n}\left(f\right)\;\longrightarrow\;1\quad
\mu-\textrm{a.e. \; uniformly\; on}\; A.\]
The set $A$ is \emph{}called \emph{uniformly} \emph{returning} if
it is uniformly returning for some $f\in P_{\mu}$.
\end{defn}
\begin{rem}
From \cite[Proposition 1.2]{KesseboehmerSlassi:05} we know that $\left(b_{n}\right)$
is regularly varying with exponent $\beta\in\left[0,1\right)$ if
and only if $\left(W_{n}\right)$ is regularly varying with the same
exponent. In this case,\[
b_{n}\sim W_{n}\Gamma\left(1-\beta\right)\Gamma\left(1+\beta\right)
\qquad\left(n\to\infty\right).\]
\end{rem}

\begin{example}\label{ex2}
Let $T:\left[0,1\right]\longrightarrow\left[0,1\right]$ be an interval
map with two increasing full branches and an indifferent fixed point
at $0$ such that there exists an absolutely continuous invariant infinite measure for $T$ like 
for instance in Example \ref{ex1}. 
In \cite{Thaler:00} Thaler introduced for this type of maps some extra conditions 
(like the convexity of $T$ in a neighbourhood of $0$ and some regularity of the density) 
which in our context guarantee that
any set $A\in\mathcal{B}_{\left[0,1\right]}$ with positive distance
from the indifferent fixed point $0$ and $\lambda\left(A\right)>0$
is uniformly returning.
\end{example}
In \cite[Proposition 1.1]{KesseboehmerSlassi:05} we have shown that any uniformly returning set
is uniform. The following proposition provides us with
conditions under which a uniform set is also uniformly returning. This will be useful 
whenever the existence 
of unifrom sets is guaranteed but Thaler's conditions in \cite{Thaler:00} are not satisfied.
This is exactly the case for the Farey map considered in 
Section \ref{sec:Application-to-continuedfraction}. 

\begin{prop}
\label{pro1}Let $A\in\mathcal{A}$ with $0<\mu\left(A\right)<\infty$
be a uniform set for $f$. If the wandering rate $\left(W_{n}\right)$
is regularly varying with exponent $1-\alpha$ for $\alpha\in\left(0,1\right]$
and the sequence $\left(\hat{T}^{n}\left(f\right)\mid_{A}\right)$
is decreasing, then $A$ is a uniformly returning set for $f$. In
this case,\[
W_{n}\hat{T}^{n}\left(f\right)\;\longrightarrow\;\frac{1}{\Gamma\left(\alpha\right)
\Gamma\left(2-\alpha\right)}\qquad\mu-\textrm{a.e. \; uniformly\; on}\; A.\]
\end{prop}
\begin{proof}
Let $\lambda,\eta\in\mathbb{R}$ be fixed but arbitrary with $0<\lambda<\eta<\infty$.
Putting\[
V_{n}:=\sum_{k=0}^{n}\hat{T}^{k}\left(f\right),\]
we have by monotonicity of $\left(\hat{T}^{n}\left(f\right)\mid_{A}\right)$
\[
\frac{\hat{T}^{\left\lfloor n\eta\right\rfloor }\left(f\right)}{V_{n}}\cdot\left(\left\lfloor n\eta\right\rfloor -\left\lfloor n\lambda\right\rfloor \right)\leq\frac{V_{\left\lfloor n\eta\right\rfloor }-V_{\left\lfloor n\lambda\right\rfloor }}{V_{n}}\leq\frac{\hat{T}^{\left\lfloor n\lambda\right\rfloor }\left(f\right)}{V_{n}}\cdot\left(\left\lfloor n\eta\right\rfloor -\left\lfloor n\lambda\right\rfloor \right).\]
Since $\left\lfloor n\eta\right\rfloor -\left\lfloor n\lambda\right\rfloor \sim n\left(\eta-\lambda\right)$
as $n\to\infty$, we have for fixed $\varepsilon\in\left(0,1\right)$
and all $n$ sufficiently large\[
n\left(1-\varepsilon\right)\left(\eta-\lambda\right)\leq\left\lfloor n\eta\right\rfloor -\left\lfloor n\lambda\right\rfloor \leq n\left(1+\varepsilon\right)\left(\eta-\lambda\right).\]
This implies for all $n$ sufficiently large\[
\frac{n\hat{T}^{\left\lfloor n\eta\right\rfloor }\left(f\right)}{V_{n}}\cdot\left(1-\varepsilon\right)\left(\eta-\lambda\right)\leq\frac{V_{\left\lfloor n\eta\right\rfloor }-V_{\left\lfloor n\lambda\right\rfloor }}{V_{n}}\leq\frac{n\hat{T}^{\left\lfloor n\lambda\right\rfloor }\left(f\right)}{V_{n}}\cdot\left(1+\varepsilon\right)\left(\eta-\lambda\right).\]
Since\[
\frac{V_{\left\lfloor n\eta\right\rfloor }-V_{\left\lfloor n\lambda\right\rfloor }}{V_{n}}\longrightarrow\eta^{\alpha}-\lambda^{\alpha}\quad\textrm{as}\; n\to\infty\quad\mu\textrm{-a.e.\quad uniformly\; on\;}A,\]
 we obtain on the one hand \[
\frac{1}{1+\varepsilon}\cdot\frac{\eta^{\alpha}-\lambda^{\alpha}}{\eta-\lambda}\leq\liminf\frac{n\hat{T}^{\left\lfloor n\lambda\right\rfloor }\left(f\right)}{V_{n}}\quad\mu\textrm{-a.e.\quad uniformly\; on\;}A.\]
Letting $\eta\to\lambda$ and $\varepsilon\to0$, it follows that\[
\alpha\lambda^{\alpha-1}\leq\liminf\frac{n\hat{T}^{\left\lfloor n\lambda\right\rfloor }\left(f\right)}{V_{n}}\quad\mu\textrm{-a.e.\quad uniformly\; on\;}A.\]
On the other hand, we obtain similarly\[
\limsup\frac{n\hat{T}^{\left\lfloor n\eta\right\rfloor }\left(f\right)}{V_{n}}\leq\alpha\eta^{\alpha-1}\quad\mu\textrm{-a.e.\quad uniformly\; on\;}A.\]
Since $\lambda\;\textrm{and}\;\eta$ are arbitrary, we have for arbitrary
$c>0$\[
\frac{n\hat{T}^{\left\lfloor nc\right\rfloor }\left(f\right)}{V_{n}}\longrightarrow\alpha c^{\alpha-1}\quad\mu\textrm{-a.e.\quad uniformly\; on\;}A.\]
Finally using $V_{\left\lfloor nc\right\rfloor }\sim c^{\alpha}V_{n}\quad\mu\textrm{-a.e.}$
uniformly on $A$ and $\left\lfloor nc\right\rfloor \sim cn$, we
obtain for $m=\left\lfloor nc\right\rfloor $\[
\frac{m\hat{T}^{m}\left(f\right)}{V_{m}}=\frac{\left\lfloor nc\right\rfloor }{n}\cdot\frac{n\hat{T}^{\left\lfloor nc\right\rfloor }\left(f\right)}{V_{n}}\cdot\frac{V_{n}}{V_{\left\lfloor nc\right\rfloor }}\longrightarrow\alpha\quad\mu\textrm{-a.e.\quad uniformly\; on\;}A.\]
From this and (\ref{eq:AsymptAaronson}) the assertion follows.
\end{proof}

\subsection{Limit laws.\label{sub:Limit-laws.}}

An important question when studying convergence in distribution for
processes defined in terms of a non-singular transformation is to
what extent the limiting behaviour depends on the initial distribution.
This is formalised as follows.

Let $\left(R_{n}\right)_{n\geq1}$ be a sequence of real valued random
variables on the $\sigma$-finite measure space $\left(X,\mathcal{A},\mu\right)$
and $R$ some random variable with values in $\left[-\infty,\infty\right]$.
Then \emph{strong distributional convergence} of $\left(R_{n}\right)_{n\geq1}$
to $R$ abbreviated by $R_{n}\stackrel{\mathcal{L\left(\mu\right)}}{\longrightarrow}R$
means that $R_{n}\stackrel{\nu}{\longrightarrow}R$ holds for all
$\nu\in\mathcal{P_{\mu}}$. In particular for $c\in\left[-\infty,\infty\right]$,\[
R_{n}\stackrel{\mathcal{L\left(\mu\right)}}{\longrightarrow}c\quad\Longleftrightarrow\quad R_{n}\stackrel{\mu}{\longrightarrow}c.\]

Now we are in the position to connect a certain renewal process with
the processes we are investigating in this paper and state the corresponding
limit laws.

Let $A\in\mathcal{A}$ be a set with $0<\mu\left(A\right)<\infty$,
and let $\left(\tau_{n}\right)_{n\in\mathbb{N}}$ be the sequence
of return times, i.e. integer valued positive random variables defined
recursively by \begin{eqnarray*}
\tau_{1}\left(x\right) & := & \varphi(x)=\inf\{ p\geq1:\; T^{p}(x)\in A\},\quad x\in X,\\
\tau_{n}\left(x\right) & := & \inf\{ p\geq1:\; T^{p+\sum_{k=1}^{n-1}\tau_{k}\left(x\right)}(x)\in A\},\quad x\in X.\end{eqnarray*}
The \emph{renewal process} is then given by\[
N_{n}(x):=\left\{ \begin{array}{ll}
\max\{ k\leq n\;:\; S_{k}\left(x\right)\leq n\}, & x\in K_{n}=\bigcup_{k=0}^{n}T^{-k}A,\\
0, & \textrm{else,}\end{array}\right.\]
where\[
S_{0}:=0,\qquad S_{n}:=\sum_{k=1}^{n}\tau_{k},\quad n\in\mathbb{N}.\]
Now we consider the so-called \emph{spent time process} $\sigma_{n}$
given by\[
\sigma_{n}(x):=\left\{ \begin{array}{ll}
n-S_{N_{n}\left(x\right)}\left(x\right), & x\in K_{n},\\
n, & \textrm{else,}\end{array}\right.\]
and the normalised \emph{spent time Kac} \emph{process}\[
\Psi_{n}:=\frac{\sum_{k=0}^{\sigma_{n}}\mu\left(A\cap\left\{ \varphi>k\right\} \right)}
{\mu\left(K_{n}\right)}=\frac{W_{\sigma_{n}}}{W_{n}}.\]
Note that\begin{equation}
S_{N_{n}\left(x\right)}\left(x\right)=Z_{n}\left(x\right):=\left\{ \begin{array}{ll}
\max\left\{ k\leq n:\; T^{k}\left(x\right)\in A\right\} , & x\in K_{n},\\
0, & \textrm{else.}\end{array}\right.\label{eq:renewalDEF}\end{equation}

We will make use of the following result from \cite{KesseboehmerSlassi:05}.
\begin{description}
\item [(UL)\textmd{~\label{theo2}~}\textmd{\emph{Uniform~law.}}]Let
$A\in\mathcal{A}$ with $0<\mu\left(A\right)<\infty$ be a uniformly
returning set. If the wandering rate $\left(W_{n}\right)$ is slowly
varying, then we have\[
\Psi_{n}\stackrel{\mathcal{L\left(\mu\right)}}{\;\longrightarrow}\;\boldsymbol{U},\]
where the random variable $\boldsymbol{U}$ is distributed uniformly
on $\left[0,1\right]$.
\end{description}
\begin{example}
Let us consider the \emph{Lasota--Yorke} map $T:\left[0,1\right]\longrightarrow\left[0,1\right]$,
defined by\[
T\left(x\right):=\left\{ \begin{array}{ll}
\frac{x}{1-x}, & x\in\left[0,\frac{1}{2}\right],\\
2x-1, & x\in\left(\frac{1}{2},1\right].\end{array}\right.\]
This map satisfies Thaler's conditions in \cite{Thaler:00}.
Hence as mentioned in Example \ref{ex2}, any compact subset $A$ of $\left(0,1\right]$ 
with $\lambda\left(A\right)>0$ is a uniformly returning set and we have\[
W_{n}\sim\log\left(n\right)\quad\textrm{as}\quad n\to\infty.\]
Consequently,\[
\frac{\log\left(\sigma_{n}\right)}{\log\left(n\right)}
\stackrel{\mathcal{L\left(\mu\right)}}{\;\longrightarrow}\;\boldsymbol{U}.\]
\end{example}

\section{Application to continued fractions\label{sec:Application-to-continuedfraction}}

In this section we will make use of the fact that one can receive
the Gauss map from the Farey map by inducing. This allows us to connect
the renewal process for the Farey map and the fluctuation process
for the continued fraction digit sum. Unfortunately the Farey map
-- unlike the Lasota--Yorke map -- does not satisfy Thaler's condition
(i)--(iv) in \cite{Thaler:00} forcing us to study this map in some detail in
order to show that the interval $(1/2,1]$ is uniformly returning.

\subsection{The Farey and Gauss map }

We consider the Farey map $T:\left[0,1\right]\rightarrow\left[0,1\right]$,
defined by\[
T\left(x\right):=\left\{ \begin{array}{ll}
T_{0}\left(x\right), & x\in\left[0,\frac{1}{2}\right],\\
T_{1}\left(x\right), & x\in\left(\frac{1}{2},1\right],\end{array}\right.\]
where\[
T_{0}\left(x\right):=\frac{x}{1-x}\qquad\textrm{and}\qquad T_{1}\left(x\right):=\frac{1}{x}-1.\]

\psfrag{A}{\large $A_{1}$}
\psfrag{gamma}{$\left(\gamma -1\right)$}
\psfrag{0}{$0$}
\psfrag{1}{$1$}
\psfrag{0.5}{$\frac{1}{2}$}
\psfrag{0.25}{$\frac{1}{4}$}
\psfrag{0.125}{$\frac{1}{5}$}
\psfrag{0.3}{$\frac{1}{3}$}
\psfrag{...}{$\cdots$}

\begin{figure} 
\centering 
\includegraphics[scale=0.45] {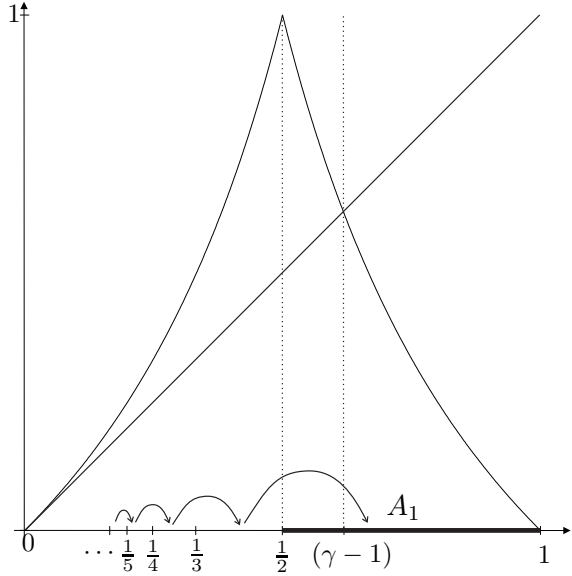} 
\caption{The Farey map $T$ and the uniformly returning set $A_{1}$.  $0$ is the critical and 
$\gamma-1$ is the non-critical fixed point of $T$, in here $\gamma$ denotes the golden ratio.
For $n\geq 2$ we have $T(1/(n+1),1/n]=(1/n    ,1/ (n-1)]$.}
\label{fig:1}  
\end{figure}
With $B\left(0\right)=\left[0,\frac{1}{2}\right]$, $B\left(1\right)=\left(\frac{1}{2},1\right]$,
and $J=\left\{ 0\right\} $, it is not difficult to see that Thaler's
condition from Example \ref{ex1} are fulfilled. It is easy to verify
that with $h\left(x\right):=\frac{d\mu}{d\lambda}\left(x\right)=\frac{1}{x}$
we have $\hat{T}\left(\1\right)=\1$ and hence $\left(\left[0,1\right],T,\mathcal{B},\mu\right)$
defines a conservative ergodic measure preserving dynamical system.
Also any Borel set $A\in\mathcal{B}$ with $\lambda\left(A\right)>0$
which is bounded away from the indifferent fixed point $0$ is a uniform
set. Furthermore, we have\[
W_{n}=\int_{\frac{1}{n+2}}^{1}\frac{1}{x}\, dx=\log\left(n+2\right)\sim\log\left(n\right)
\qquad\left(n\to\infty\right).\]
The inverse branches of the Farey map are \begin{eqnarray*}
u_{0}\left(x\right) & := & \left(T_{0}\right)^{-1}\left(x\right)={\displaystyle \frac{x}{1+x}},\\
u_{1}\left(x\right) & := & \left(T_{1}\right)^{-1}\left(x\right)={\displaystyle \frac{1}{1+x}}.\end{eqnarray*}
For $x\neq0$ the map $u_{0}\left(x\right)$ is conjugated to the
right translation $x\mapsto F\left(x\right):=x+1$, i.e.\[
u_{0}=J\circ F\circ J\qquad\textrm{with}\quad J\left(x\right)=J^{-1}\left(x\right)=\frac{1}{x}.\]
This shows that for the $n$-th iterate we have\begin{equation}
u_{0}^{n}\left(x\right)=J\circ F^{n}\circ J\left(x\right)=\frac{x}{1+nx}.\label{u0}\end{equation}
Moreover, we have $u_{1}\left(x\right)=J\circ F\left(x\right)$. 

Let $\mathcal{F}=\left\{ A_{n}\right\} _{n\geq1}$ be the countable
collection of pairwise disjoint subintervals of $\left[0,1\right]$
given by $A_{n}=\left(\frac{1}{n+1},\frac{1}{n}\right]$. Setting
$A_{0}=\left[0,1\right)$, it is easy to check that $T\left(A_{n}\right)=A_{n-1}$
for all $n\geq1$. The \emph{first entry time} $e:\I\rightarrow\mathbb{N}$
in the interval $A_{1}$ is defined as\[
e\left(x\right):=\min\left\{ k\geq0:\; T^{k}\left(x\right)\in A_{1}\right\} .\]
Then the first entry time is connected to the first digit in the continued
fraction expansion by\[
a_{1}\left(x\right)=1+e\left(x\right)\quad\textrm{and}\quad\varphi\left(x\right)=a_{1}\circ T\left(x\right),\quad x\in\I.\]
We now consider the \emph{induced map} $S:\I\rightarrow\I$ defined
by\[
S\left(x\right):=T^{e\left(x\right)+1}\left(x\right).\]
Since for all $n\ge1$ \[
\left\{ x\in\I:\; e\left(x\right)=n-1\right\} =A_{n}\cap\I,\]
we have by (\ref{u0}) for any $x\in A_{n}\cap\I$ \[
S\left(x\right)=T^{n}\left(x\right)=T_{1}\circ T_{0}^{n-1}\left(x\right)=\frac{1}{x}-n=\frac{1}{x}-a_{1}(x).\]
This implies that the induced transformation $S$ coincides with Gauss
map $G$ on $\I$. 

In the next lemma we connect the number theoretical process $X_{n}$
defined in (\ref{process}) with the renewal process $Z_{n}$ with
respect to the Farey map defined in (\ref{eq:renewalDEF}).
\begin{lem}
\label{lemma1}Let $A_{1}:=\left(\frac{1}{2},1\right]$ and $K_{n}:=\bigcup_{k=0}^{n}T^{-k}A_{1}$.
Then for the process $X_{n}$ defined in (\ref{process}) we have
for all $x\in\I$ \[
X_{n}\left(x\right)=\left\{ \begin{array}{ll}
1+Z_{n-1},\qquad & x\in K_{n-1},\\
0, & \textrm{else.}\end{array}\right.\]
\end{lem}
\begin{proof}
From the above discussion and the definition of the renewal theoretic
process $S_{N_{n}}$ in Subsection \ref{sub:Limit-laws.} we deduce the following.
\begin{itemize}
\item For $x\in\I\cap K_{n-1}^{C}$ we have $a_{1}\left(x\right)>n$, implying
$X_{n}\left(x\right)=0$.
\item For $x\in\I\cap K_{n-1}$ we distinguish two cases. Either the process
starts in $x\in A_{1}$, then we have $a_{1}\left(x\right)=1$ and
inductively for $n\ge2$ \[
a_{n}\left(x\right)=\tau_{n-1}\left(x\right),\]
 or the process starts in $x\in A_{1}^{C}$, then we have $a_{1}\left(x\right)=
1+\tau_{1}\left(x\right)$
and inductively for $n\ge2$\[
a_{n}\left(x\right)=\tau_{n}\left(x\right).\]
\end{itemize}
From this the assertion follows. 
\end{proof}
To show that $A_{1}$ is uniformly returning we need the following
lemma.
\begin{lem}
\label{lem1}Let\[
\mathcal{D}:=\left\{ f\in\mathcal{P}_{\mu}:\; f\in\mathcal{C}^{2}\left(\left(0,1\right)\right)\;\textrm{with}\; f'>0\;\textrm{and}\; f''\leq0\right\} .\]
Then we have for all $n\in\mathbb{N}$\[
f\in\mathcal{D}\quad\implies\quad\hat{T}^{n}\left(f\right)\in\mathcal{D}.\]
\end{lem}
\begin{proof}
First note, that it suffices to show the assertion only for $n=1$
since the lemma then follows by induction. For $f\in\mathcal{D}$
we have\[
\hat{T}\left(f\right)=\frac{1}{h}\cdot\hat{P}\left(h\cdot f\right),\]
 where $\hat{P}$ denotes the Perron-Frobenius operator of $T$ restricted
to $L_{1}\left(\lambda\right)$. Using the inverse branches of $T$
this operator is given by\[
\hat{P}\left(g\right)=g\circ u_{0}\cdot\left|u_{0}'\right|+g\circ u_{1}\cdot\left|u_{1}'\right|,\quad\textrm{for\; all}\; g\in L_{1}\left(\lambda\right).\]
It follows that\[
\hat{T}\left(f\right)\left(x\right)=\frac{f\left(\frac{x}{x+1}\right)+xf\left(\frac{1}{x+1}\right)}{x+1},\quad\textrm{for\; all}\quad x\in\left[0,1\right].\]
Hence, also $\hat{T}\left(f\right)$ is differentiable on $\left(0,1\right)$
and by the monotonicity of $f$ and $f'$ we have\[
\hat{T}\left(f\right)'\left(x\right)=\underbrace{\frac{f'\left(\frac{x}{x+1}\right)
-xf'\left(\frac{1}{x+1}\right)}{\left(x+1\right)^{3}}}_{>0}+
\underbrace{\frac{f\left(\frac{1}{x+1}\right)-f\left(\frac{x}{x+1}\right)}
{\left(x+1\right)^{2}}}_{>0}\]
implying $\hat{T}\left(f\right)'>0$. 
Furthermore, an easy calculation shows\begin{eqnarray*}
\hat{T}\left(f\right)''\left(x\right) & = & \frac{f''\left(\frac{x}{x+1}\right)
+xf''\left(\frac{1}{x+1}\right)}{\left(x+1\right)^{5}}+
\frac{2\left(f\left(\frac{x}{x+1}\right)-f\left(\frac{1}{x+1}\right)\right)}
{\left(x+1\right)^{3}}\\
 &  & \qquad+\frac{2\left(x+1\right)\left(\left(x-1\right)f'\left(\frac{1}{x+1}\right)
-2f'\left(\frac{x}{x+1}\right)\right)}{\left(x+1\right)^{5}}\leq0.\end{eqnarray*}
This finishes the proof. 
\end{proof}
We remark that $\mathcal{D}\neq\emptyset$, since $\textrm{id}_{\left[1,0\right]}\in\mathcal{D}$.
\begin{lem}
\label{lemma2}The set $A_{1}=\left(\frac{1}{2},1\right]$ is uniformly
returning for any $f\in\mathcal{D}$.
\end{lem}
\begin{proof}
Since $A_{1}$ is uniform for any $f\in\mathcal{P}_{\mu}$ (cf. Example
\ref{ex1}) we have in view of Proposition \ref{pro1} only to verify
that $\left(\hat{T}^{n}\left(f\right)\mid_{A_{1}}\right)$ is decreasing.
For all $x\in A_{1}$, $n\in\mathbb{N}_{0}$, we have\begin{eqnarray*}
\hat{T}^{n+1}\left(f\right)\left(x\right) & = & \frac{\hat{T}^{n}\left(f\right)\left(\frac{x}{x+1}\right)+x\hat{T}^{n}\left(f\right)\left(\frac{1}{x+1}\right)}{x+1}\\
 & = & \frac{1}{x+1}\hat{T}^{n}\left(f\right)\left(\frac{x}{x+1}\right)+\frac{x}{x+1}\hat{T}^{n}\left(f\right)\left(\frac{1}{x+1}\right)\end{eqnarray*}
Since by Lemma \ref{lem1} each $\hat{T}^{n}\left(f\right)$ is concave
and increasing on $\left[0,1\right]$ and $\frac{1}{x+1}+\frac{x}{x+1}=1$,
we have for all $x\geq\sqrt{2}-1$,\[
\hat{T}^{n+1}\left(f\right)\left(x\right)\leq\hat{T}^{n}\left(f\right)\left(\frac{2x}
{\left(x+1\right)^{2}}\right)\leq\hat{T}^{n}\left(f\right)\left(x\right).\]
 This finishes the proof.
\end{proof}

\subsection{Distributional limit laws for digit sums\label{sub:Distributional-limit-laws}}

With the above preparations we can now state the main results.
\begin{thm}
\label{thm:main}Let $X_{n}$ be the process given in (\ref{process}).
Then the following holds.
\begin{enumerate}
\item We have \begin{equation}
\frac{\log\left(n-X_{n}\right)}{\log\left(n\right)}\stackrel{\mathcal{L\left(\mu\right)}}
{\longrightarrow}\boldsymbol{U},\label{gl3}\end{equation}
 where the random variable $\boldsymbol{U}$ is distributed uniformly
on $\left[0,1\right]$. 
\item For $f\in\mathcal{D}$ set $d\nu:=f\,d\mu$. Then for any $a\in\left(0,1\right)$ we have 
\begin{equation}
\nu\left(\frac{n-X_{n}}{n}>a\right)\sim\frac{-\log\left(a\right)}
{\log\left(n\right)}\quad\textrm{as}\quad n\to\infty.
\label{gl2}\end{equation}
\end{enumerate}
\end{thm}
\begin{proof}
Using Lemma \ref{lemma2} and the fact that $W_{n}\sim\log\left(n\right)$
we have by (UL) \[
\frac{W_{n-Z_{n}}}{\log n}\stackrel{\mathcal{L\left(\mu\right)}}{\longrightarrow}\boldsymbol{U}.\]
Consequently $n-Z_{n}\stackrel{\mathcal{L\left(\mu\right)}}{\longrightarrow}\infty$
and hence \[
\frac{W_{n-Z_{n}}}{\log\left(n-Z_{n}\right)}\stackrel{\mu}{\longrightarrow}1.\]
Thus by Lemma \ref{lemma1}, the convergence in (\ref{gl3}) holds. 

For the second part of the theorem let $f\in\mathcal{D}$ and $a\in\left(0,1\right)$
be fixed and set  $d\nu:=f\,d\mu$. It is not difficult to verify that\[
\nu\left(\frac{n-Z_{n}}{n}>a\right)\sim\nu\left(\frac{n-X_{n}}{n}>a\right)\quad\textrm{as}\quad n\to\infty.\]
Therefore, to prove (\ref{gl2}) it suffices to show\[
\nu\left(\frac{n-Z_{n}}{n}>a\right)\sim\frac{-\log\left(a\right)}
{\log\left(n\right)}\quad\textrm{as}\quad n\to\infty.\]
In fact, we have\begin{eqnarray*}
\nu\left(\frac{n-Z_{n}}{n}>a\right) 
& = & \sum_{k=0}^{\left\lfloor n\left(1-a\right)\right\rfloor }\nu
\left(K_{n}\cap\left\{ Z_{n}=k\right\} \right)\\
& = & \sum_{k=0}^{\left\lfloor n\left(1-a\right)\right\rfloor }\nu
\left(T^{-k}\left(A_{1}
\cap\left\{ \varphi>n-k\right\} \right)\right)\\
& = & \sum_{k=0}^{\left\lfloor n\left(1-a\right)\right\rfloor }\int_{A_{1}}\1_{A_{1}
\cap\left\{ \varphi>n-k\right\} }\hat{T}^{k}\left(f\right)\; d\mu.
\end{eqnarray*}
Let $\delta\in\left(0,1-a\right)$
and $\varepsilon\in\left(0,1\right)$ be fixed but arbitrary and divide
the above sum into two parts as follows.\[
\nu\left(\frac{n-Z_{n}}{n}>a\right)=\sum_{k=0}^{\left\lfloor n\delta\right\rfloor -1}\cdots+\sum_{k=\left\lfloor n\delta\right\rfloor }^{\left\lfloor n\left(1-a\right)\right\rfloor }\cdots=:\; I\left(n\right)+J\left(n\right).\]
We first note, that since $A_{1}\cap\left\{ \varphi>n\right\} =\left[\frac{n+2}{n+3},1\right]$
we have \begin{equation}
\mu\left(A_{1}\cap\left\{ \varphi>n\right\} \right)=\int_{\frac{n+2}{n+3}}^{1}\frac{1}{x}\: dx\;\sim\;\frac{1}{n}\quad\textrm{as}\quad n\to\infty.\label{u1}\end{equation}
Also, by monotonicity of $\left(\1_{A_{1}\cap\left\{ \varphi>n\right\} }\right)$
we have\[
I\left(n\right)\leq\int_{A_{1}}\1_{A_{1}\cap\left\{ \varphi>n+1-\left\lfloor n\delta\right\rfloor \right\} }\sum_{k=0}^{\left\lfloor n\delta\right\rfloor -1}\hat{T}^{k}\left(f\right)\; d\mu.\]
 Combining both observation and using the fact that $A_{1}$ is uniform
for $f$ (cf. Example \ref{ex1}) and that (\ref{eq:AsymptAaronson})
holds, we obtain for sufficiently large $n$\begin{eqnarray*}
I\left(n\right) & \leq & \left(1+\varepsilon\right)^{2}\frac{\left\lfloor n\delta\right\rfloor -1}{n-\left\lfloor n\delta\right\rfloor +1}\cdot\frac{1}{\log\left(\left\lfloor n\delta\right\rfloor -1\right)}\\
 &  & \sim\left(1+\varepsilon\right)^{2}\frac{\delta}{1-\delta}\cdot\frac{1}{\log\left(n\right)}\quad\textrm{as}\quad n\to\infty.\end{eqnarray*}
Thus,\[
\limsup_{n\to\infty}\log\left(n\right)\cdot I\left(n\right)\leq\left(1+\varepsilon\right)^{3}\frac{\delta}{1-\delta}.\]
Letting $\delta\to0$, we conclude\begin{equation}
I\left(n\right)=o\left(\frac{1}{\log\left(n\right)}\right),\quad n\to\infty.\label{log1}\end{equation}
For the second part of the sum we have to show that \begin{equation}
J\left(n\right)\sim\frac{-\log\left(a\right)}{\log\left(n\right)}\quad\textrm{as}\quad n\to\infty.\label{j1}\end{equation}
A similarly argument as in the proof of  Lemma 3.3 in \cite{KesseboehmerSlassi:05}
shows that for all $n$ sufficiently large and $k\in\left[\left\lfloor n\delta\right\rfloor ,\left\lfloor n\left(1-a\right)\right\rfloor \right]$
we have uniformly on $A_{1}$\begin{equation}
\left(1+\varepsilon\right)\frac{1}{\log\left(n\right)}\leq\hat{T}^{k}\left(f\right)\leq\left(1+\varepsilon\right)^{2}\frac{1}{\log\left(n\right)}.\label{f1}\end{equation}
Hence, using the right-hand side of (\ref{f1}) and the asymptotic
formula (\ref{u1}), we obtain for $n$ sufficiently large \begin{eqnarray*}
J\left(n\right) & \leq & \frac{\left(1+\varepsilon\right)^{2}}{\log\left(n\right)}\cdot\sum_{k=n-\left\lfloor n\left(1-a\right)\right\rfloor }^{n-\left\lfloor n\delta\right\rfloor }\mu\left(A_{1}\cap\left\{ \varphi>k\right\} \right)\\
 &  & \sim\frac{\left(1+\varepsilon\right)^{3}}{\log\left(n\right)}\cdot\log\left(\frac{1-\delta}{a}\right)\quad\textrm{as}\quad n\to\infty,\end{eqnarray*}
which implies \[
\limsup_{n\to\infty}\log\left(n\right)\cdot J\left(n\right)\leq\left(1+\varepsilon\right)^{4}\log\left(\frac{1-\delta}{a}\right).\]
Similarly, using the left-hand side of (\ref{f1}), we get\[
\liminf_{n\to\infty}\log\left(n\right)\cdot J\left(n\right)\geq\left(1+\varepsilon\right)^{3}\log\left(\frac{1-\delta}{a}\right).\]
Since $\varepsilon$ and $\delta$ were arbitrary, (\ref{j1}) holds. 
Combining (\ref{log1}) and (\ref{j1}) proves the second part of
the theorem.
\end{proof}
Finally, since id$_{[0,1]}\in\mathcal{D}$, Theorem \ref{thm:main} in
particular implies the claims for $x\cdot d\mu(x)=d\lambda(x)$ as stated in the introduction. Also note
that since $\lambda\sim\mu$, by (\ref{eq:Negelg}) we have\[
\frac{n-X_{n}}{n}\stackrel{\mathcal{L\left(\mu\right)}}{\longrightarrow}0,\]
 and it follows directly from (\ref{gl3}) that\[
n-X_{n}\stackrel{\mathcal{L\left(\mu\right)}}{\longrightarrow}\infty.\]

%\bibliographystyle{alpha}
%\bibliography{InfinitErgod,InfinitErgodic,InfiniteErgod}

\end{document}